\documentclass[amstex,12pt,russian,amssymb]{article}

\usepackage{mathtext}
\usepackage[cp1251]{inputenc}
\usepackage[T2A]{fontenc}
\usepackage[russian]{babel}
\usepackage[dvips]{graphicx}
\usepackage{amsmath}
\usepackage{amssymb}
\usepackage{amsxtra}
\usepackage{latexsym}
\usepackage{ifthen}

\begin{document}

\newcounter{lemma}
\newcommand{\lemma}{\par \refstepcounter{lemma}%
{\bf Лемма \arabic{lemma}.}}

\newcounter{corollary}
\newcommand{\corollary}{\par \refstepcounter{corollary}%
{\bf Следствие \arabic{corollary}.}}

\newcounter{remark}
\newcommand{\remark}{\par \refstepcounter{remark}%
{\bf Замечание \arabic{remark}.}}

\newcounter{theorem}
\newcommand{\theorem}{\par \refstepcounter{theorem}%
{\bf Теорема \arabic{theorem}.}}

\newcounter{proposition}
\newcommand{\proposition}{\par \refstepcounter{proposition}%
{\bf Предложение \arabic{proposition}.}}

\renewcommand{\refname}{\centerline{\bf Список литературы}}

\newcommand{\proof}{{\it Доказательство.\,\,}}

\noindent УДК 517.5

{\bf Е.А.~Севостьянов} (Житомирский государственный университет им.\
И.~Франко)

\medskip
{\bf Є.О.~Севостьянов} (Житомирський державний університет ім.\
І.~Франко)

\medskip
{\bf E.A.~Sevost'yanov} (Ivan Franko Zhitomir State University)

\medskip
{\bf О равностепенной непрерывности гомеоморфизмов в замыкании
области метрического пространства}

\medskip
{\bf Про одностайну неперервність гомеоморфізмів в замиканні області
метричного простору}

\medskip
{\bf On equicontinuity of homeomorphisms in a closure of a domain in
metric spaces}

\medskip
В работе изучаются гомеоморфизмы на метрических пространствах,
обобщающие квазиконформные отображения. Установлено, что семейства
указанных отображений при определённых условиях на метрические
пространства и границы соответствующих областей являются
равностепенно непрерывными в замыкании заданной области.

\medskip
В роботі вивчаються гомеоморфізми на метричних просторах, що
узагальнюють квазіконформні відображення.  Встановлено, що сім'ї
вказаних відображень за деяких умов на метричні простори і межі
відповідних областей є одностайно неперервними в замиканні заданої
області.

\medskip
In the present paper, homeomorphisms in metric spaces which
generalize quasiconformal mappings, are investigated. It is proved
that, at some conditions on metric spaces and boundaries of
corresponding domains, families of above mappings are equicontinuous
in a closure of a given domain.

\newpage
{\bf 1. Введение.} В настоящей работе речь идёт об изучении
гомеоморфизмов, представляющих собой квазиконформные отображения с
весом (см. \cite{AC$_1$}, \cite{BGMV}, \cite{Cr$_1$}, \cite{GG},
\cite{GRY} и \cite{RSa}). Здесь мы исследуем локальное поведение так
называемых кольцевых $Q$-го\-ме\-о\-мор\-физ\-мов в области
метрического пространства, а также их поведение в замыкании области
метрического пространства (см., напр., \cite{RS} и \cite{Sev$_1$}).
Основоположниками указанных исследований являются Някки и Палка,
полностью описавшие ситуацию квазиконформных отображений в
евклидовом случае (\cite{NP}, \cite{Na}). Для отображений с
неограниченной характеристикой аналогичные вопросы рассмотрены
автором в евклидовом пространстве \cite{Sev$_2$}. Более того,
частично рассмотрена и ситуация метрических пространств (см.
\cite{Sev$_4$}), где некоторые из приведённых ниже утверждений
установлены в частном случае. Отметим, кроме того, статью
\cite{RSa}, где получены значительные результаты о граничном
поведении $Q$-го\-ме\-о\-мор\-физ\-мов. Хотя здесь не изучались
вопросы, связанные с локальным и глобальным поведением отображений,
именно здесь данный объект ($Q$-гомеоморфизмы) впервые рассмотрен в
метрических пространствах.

\medskip
Перейдём к изложению содержательной части работы и её основных
результатов. Хорошо известно, что отображения классов Соболева и
Орлича--Соболева в пространстве ${\Bbb R}^n$ удовлетворяют
соотношениям вида
\begin{equation}\label{eq6}
M(f(\Gamma))\leqslant
\int\limits_{\varepsilon<|x-x_0|<\varepsilon_0}Q(x)\eta^n(|x-x_0|)dm(x)
\end{equation}
для произвольной измеримой по Лебегу функции $\eta:[\varepsilon,
\varepsilon_0]\rightarrow [0, \infty]$ такой, что
$\int\limits_{\varepsilon}^{\varepsilon_0}\eta(t)dt\geqslant 1,$ где
$\Gamma$ -- семейство кривых, соединяющих сферы с центром в точке
$x_0$ и радиусов $\varepsilon$ и $\varepsilon_0,$ $m$ -- мера Лебега
в ${\Bbb R}^n,$ а $M$ -- конформный модуль семейств кривых (см.,
напр., \cite[следствие~5]{KRSS} и \cite[теорема~1]{Sev$_3$}). Наша
ближайшая цель -- изучить некоторые свойства отображений аналогичным
тем, что удовлетворяют соотношениям (\ref{eq6}), в метрических
пространствах. В настоящей статье основное внимание уделяется
локальному поведению отображений и поведению их в замыкании заданной
области.

\medskip
В дальнейшем $(X, d, \mu)$ и $\left(X^{\,\prime}, d^{\,\prime},
\mu^{\,\prime}\right)$ -- метрические пространства с метриками $d$ и
$d^{\,\prime}$ и борелевскими мерами $\mu$ и $\mu^{\,\prime}.$ Под
{\it областью $G$} в пространстве $X$ следует понимать {\it открытое
линейно связное множество} в $X.$ {\it Кривой} $\gamma$ в $X$
называется непрерывное отображение $\gamma:[a, b]\rightarrow X.$
{\it Длиной кривой $\gamma$ на отрезке $[a, b]$ называется величина}
$$l(\gamma):=\sup\sum\limits_{i=1}^n d(\gamma(t_i),\gamma(t_{i-1}))\,,$$
где $\sup$ берётся по всем возможным разбиениям $a=t_0\leqslant
t_1\leqslant\ldots\leqslant t_n:=b.$ Если $l(\gamma)<\infty,$ кривая
называется спрямляемой и, значит, корректно определена функция длины
$s_{\gamma}(t),$ означающая длину кривой $\gamma_{[a, t]},$ $t\in
[a, b].$ В таком случае, имеет место представление
$$\gamma(t)=\gamma^0\circ s_{\gamma}(t)\,,$$
где $\gamma^0$ называется {\it нормальным представлением} кривой
$\gamma$ (см., напр., \cite[гл.7, соотношение~(7.2)]{He}). {\it
Интегралом} от борелевской функции $\rho:G\rightarrow [0, \infty]$
называется величина
$$\int\limits_{\gamma}\rho(x)|dx|=\int\limits_{0}^{l(\gamma)}\rho(\gamma^0(t))dt\,.$$
Под семейством кривых $\Gamma$ подразумевается некоторый
фиксированный набор кривых $\gamma.$ Борелева функция
$\rho:X\,\rightarrow [0,\infty]$ называется {\it допустимой} для
семейства $\Gamma$ кривых $\gamma$ в $X,$ если
%
%
$$\int\limits_{\gamma}\rho (x)|dx|\geqslant 1$$
%
%
для всех (локально спрямляемых) кривых $ \gamma \in \Gamma$ (т.е.,
произвольная кривая $\gamma$ семейства $\Gamma$ имеет длину, не
меньшую $1$ в метрике $\rho$). В этом случае мы пишем: $\rho \in
{\rm adm} \,\Gamma.$

\medskip Пусть $p\geqslant 1$ -- фиксированное число, тогда
{\it $p$-модулем} семейства кривых $\Gamma $ называется величина
$$M_p(\Gamma)=\inf_{\rho \in \,{\rm adm}\,\Gamma}
\int\limits_X \rho ^p (x) d\mu(x)\,.$$
При этом, если ${\rm adm}\,\Gamma=\varnothing,$ то полагаем:
$M_{p}(\Gamma)=\infty$ (см. \cite[разд.~6 на с.~16]{Va}).

\medskip
Говорят, что семейство кривых $\Gamma_1$ {\it минорируется}
семейством $\Gamma_2,$ пишем $\Gamma_1\,>\,\Gamma_2,$
если для каждой кривой $\gamma\,\in\,\Gamma_1$ существует подкривая,
которая принадлежит семейству $\Gamma_2.$
В этом случае,
\begin{equation}\label{eq32*A}
\Gamma_1
> \Gamma_2 \quad \Rightarrow \quad M_p(\Gamma_1)\leqslant M_p(\Gamma_2)
\end{equation} (см. \cite[теорема~1]{Fu}).

\medskip
Пусть $E,$ $F\subset X$ -- произвольные множества. В дальнейшем
через $\Gamma(E,F,X)$ мы обозначаем семейство всех кривых
$\gamma:[a,b]\rightarrow X,$ которые соединяют $E$ и $F$ в $X,$ т.е.
$\gamma(a)\in E,\,\gamma(b)\in F$ и $\gamma(t)\in X$ при
$t\in(a,\,b).$ Пусть $G$ и $G^{\,\prime}$ -- области с конечными
хаусдорфовыми размерностями $\alpha\geqslant 2$ и
$\alpha^{\,\prime}\geqslant 2$ в метрических пространствах $(X, d,
\mu)$ и $(X^{\,\prime}, d^{\,\prime}, \mu^{\,\prime}),$
соответственно, и пусть $Q:G\rightarrow [0, \infty]$ -- измеримая
функция. Всюду далее
$$B(x_0, r)=\left\{x\in X: d(x,
x_0)<r\right\},
S(x_0, r)=\left\{x\in X: d(x, x_0)=r\right\}\,,$$
\begin{equation}\label{eq49***}
A=A(x_0, r_1, r_2)=\left\{x\in X: r_1<d(x, x_0)<r_2\right\}\,.
\end{equation}
Зафиксируем $p, q\geqslant 1$ и будем называть гомеоморфизм
$f:G\rightarrow G^{\,\prime}$  {\it кольцевым $Q$-гомеоморфизмом в
точке $x_0\in G$ относительно $p$ и $q$-модулей}, если при любых
$0<r_1<r_2<{\rm dist}\, (x_0,
\partial G)$ и для любых сфер $S_1=S(x_0, r_1)$ и $S_2=S(x_0, r_2)$
выполнено неравенство
\begin{equation}\label{eq1}
M_p(f(\Gamma(S_1, S_2, A)))\leqslant \int\limits_A Q(x)\cdot
\eta^q(d(x, x_0))d\mu(x)
\end{equation}
для каждой измеримой функции $\eta \colon  (r_1,r_2)\rightarrow
[0,\infty ]\,$ такой, что
\begin{equation}\label{eq*3!!}
\int\limits_{r_1}^{r_2}\eta(r)dr\geqslant 1\,,\end{equation}
где $A=A(x_0, r_1, r_2).$

\medskip
Пусть $\left(X,\,d\right)$ и
$\left(X^{\,{\prime}},{d}^{\,{\prime}}\right)$ --- метрические
пространства с расстояниями  $d$  и ${d}^{\,{\prime}},$
соответственно. Семейство $\frak{F}$ непрерывных отображений
$f\colon X\rightarrow {X}^{\,\,\prime}$ называется {\it нормальным},
если из любой последовательности отображений $f_{m} \in \frak{F}$
можно выделить подпоследовательность $f_{m_{k}}$, которая сходится
локально равномерно в $X$ (т.е., равномерно на любых компактных
подмножествах $X$) к непрерывной функции $f\colon \,X\,\rightarrow\,
X^{\,\prime}.$

\medskip
Введенное  понятие  очень  тесно  связано  со  следующим. Семейство
$\frak{F}$ отображений $f\colon X\rightarrow {X}^{\,\prime}$
называется {\it равностепенно непрерывным в точке} $x_0 \in X,$ если
для любого $\varepsilon>0$ найдётся такое $\delta>0$, что
${d}^{\,\prime} \left(f(x),f(x_0)\right)<\varepsilon$ для всех таких
$x,$ что $d(x,x_0)<\delta$ и для всех $f\in \frak{F}.$ Говорят, что
$\frak{F}$ {\it равностепенно непрерывно}, если $\frak{F}$
равностепенно непрерывно в каждой  точке $x_0\in X.$ Согласно одной
из версий теоремы Арцела--Асколи (см., напр.,
\cite[пункт~20.4]{Va}), если $\left(X,\,d\right)$ --- сепарабельное
метрическое пространство, а $\left(X^{\,\prime},\,
d^{\,\prime}\right)$ --- компактное метрическое пространство, то
семейство $\frak{F}$ отображений $f\colon X\rightarrow
{X}^{\,\prime}$ нормально тогда  и только тогда, когда  $\frak{F}$
равностепенно непрерывно.

\medskip
Следующее определение может быть найдено, напр., в
\cite[разд.~4]{RSa}. Будем говорить, что интегрируемая в $B(x_0, r)$
функция ${\varphi}\colon G\rightarrow{\Bbb R}$ имеет {\it конечное
среднее колебание} в точке $x_0\in G$, $\varphi\in FMO(x_0),$ если
%
%
%
%
$$\limsup\limits_{\varepsilon\rightarrow 0}\frac{1}{\mu(B(
x_0,\,\varepsilon))}\int\limits_{B(x_0,\,\varepsilon)}
|{\varphi}(x)-\overline{{\varphi}}_{\varepsilon}|\,
d\mu(x)<\infty,$$
%
%
где
$\overline{{\varphi}}_{\varepsilon}=\frac{1} {\mu(B(
x_0,\,\varepsilon))}\int\limits_{B(x_0,\,\varepsilon)}
{\varphi}(x)\, d\mu(x).$

\medskip
Пусть $(X, d, \mu)$ --- метрическое пространство с метрикой $d,$
наделённое локально конечной борелевской мерой $\mu.$ Следуя
\cite[раздел 7.22]{He} будем говорить, что борелева функция
$\rho\colon  X\rightarrow [0, \infty]$ является {\it верхним
градиентом} функции $u\colon X\rightarrow {\Bbb R},$ если для всех
спрямляемых кривых $\gamma,$ соединяющих точки $x$ и $y\in X$
выполняется неравенство $|u(x)-u(y)|\leqslant
\int\limits_{\gamma}\rho\,|dx|.$ Будем также говорить, что в
указанном пространстве $X$ выполняется $(1; p)$-неравенство
Пуанкаре, если найдутся постоянные $C\geqslant 1$ и $\tau\geqslant
1$ так, что для каждого шара $B\subset X,$ произвольной ограниченной
непрерывной функции $u\colon \tau B\rightarrow {\Bbb R}$ и любого её
верхнего градиента $\rho$ выполняется следующее неравенство:
$$\frac{1}{\mu(B)}\int\limits_{B}|u-u_B|d\mu(x)\leqslant C\cdot({\rm diam\,}B)\left(\frac{1}{\mu(\tau B)}
\int\limits_{\tau B}\rho^p d\mu(x)\right)^{1/p}\,,$$
где $u_B:=\frac{1}{\mu(B)}\int\limits_{B}u d\mu(x).$ Метрическое
пространство $(X, d, \mu)$ назовём {\it $\widetilde{Q}$-регулярным
по Альфорсу} при некотором $\widetilde{Q}\geqslant 1,$ если при
каждом $x_0\in X,$ некоторой постоянной $C\geqslant 1$ и
произвольного $R<{\rm diam}\,X$
$$\frac{1}{C}R^{\widetilde{Q}}\leqslant \mu(B(x_0, R))\leqslant CR^{\widetilde{Q}}.$$
Здесь иногда берутся замкнутые шары $\overline{B(x_0, R)},$ что
ввиду предельных свойств меры не является принципиальным (см.
\cite[теорема~9.1, гл.~I]{Sa}). Как известно, $\alpha$-регулярные по
Альфорсу пространства имеют хаусдорфову размерность $\alpha$ (см.
\cite[с.~61--62]{He}). Более того, нетрудно видеть, что в таких
пространствах области $G$ также имеют хаусдорфову размерность
$\alpha$ (см. там же). Условимся говорить, что метрическое
пространство $X$ {\it локально связно}, если для произвольной
окрестности $U$ произвольной точки $x_0\in X$ найдётся окрестность
$V\subset U,$ являющаяся связной (см. \cite[I.49.6]{Ku}).
Справедлива следующая

\medskip
\begin{theorem}\label{theor4*!} {\sl\,Пусть $G$ -- область в локально связном и локально компактном
метрическом пространстве $(X, d, \mu)$ с конечной хаусдорфовой
размерностью $\alpha\geqslant 2,$ а $(X^{\,\prime}, d^{\,\prime},
\mu^{\,\prime})$ -- метрическое пространство, которое является
$\alpha^{\,\prime}$-регулярным по Альфорсу, и в котором выполнено
$(1; p)$-неравенство Пуанкаре, $p\in (\alpha^{\,\prime}-1,
\alpha^{\,\prime}].$

Пусть $B_R\subset X^{\,\prime}$
--- некоторый фиксированный шар радиуса $R.$ Обозначим
через $\frak{R}_{p, q, x_0, Q, B_R, \delta}(G)$ семейство кольцевых
$Q$-го\-ме\-о\-мор\-физ\-мов $f\colon G\rightarrow B_R\setminus K_f$
в точке $x_0\in G$ относительно $p$ и $q$-модулей области $G$ на
некоторую область $f(G)\subset X^{\,\prime}$ таких, что
$$\sup\limits_{x, y\in K_f} d^{\,\prime}(x, y)\geqslant
\delta>0\,,$$ где $K_f\subset B_R$ -- некоторый континуум. Тогда
семейство отображений $\frak{R}_{p, q, x_0, Q, B_R, \delta}(G)$
является равностепенно непрерывным в точке $x_0\in G,$ если $q\in
(1, \alpha]$ и $Q\in FMO(x_0).$}
\end{theorem}

\medskip
Как будет видно из дальнейших рассуждений, приведённый результат о
равностепенной непрерывности отображений распространяется также на
точки замыкания области. При этом, здесь требуются определённые
дополнительные условия на границу области. В связи с этим, напомним
некоторые определения. Область $D$ называется {\it локально связной
в точ\-ке} $x_0\in\partial D,$ если для любой окрестности $U$ точки
$x_0$ найдется окрестность $V\subset U$ такая, что множество $V\cap
D$ связно (см. \cite[I.49.6]{Ku}). Аналогично, область $D$ будет
называться {\it локально линейно связной в точ\-ке} $x_0\in\partial
D,$ если для любой окрестности $U$ точки $x_0$ найдется окрестность
$V\subset U$ такая, что множество $V\cap D$ линейно связно. Согласно
\cite{GM}, область $D$ в ${\Bbb R}^n$ будем называть {\it областью
квазиэкстремальной длины относительно $p$-модуля}, сокр. {\it
$QED$-облас\-тью относительно $p$-модуля}, если
\begin{equation}\label{eq4***}
M_p(\Gamma(E, F, X))\leqslant A\cdot M_p(\Gamma(E, F, D))
\end{equation}
для конечного числа $A\geqslant 1$ и всех континуумов $E$ и $F$ в
$D.$ Пусть $p,q\geqslant 1,$ $D$ -- область в метрическом
пространстве $(X, d, \mu)$ с конечной хаусдорфовой размерностью
$\alpha\geqslant 2,$ $Q:X\rightarrow [0, \infty]$ -- измеримая
функция, $Q(x)\equiv 0$ при всех $x\not\in D.$ Отображение
$f:D\rightarrow X$ будем называть {\it кольцевым
$Q$-отоб\-ра\-же\-нием в точке $x_0\in
\partial D$ относительно $p$ и $q$-модулей,} если для некоторого $r_0=r(x_0)$, произвольных
кольца $A=A(x_0,r_1,r_2),$ центрированного в точке $x_0,$ радиусов:
$r_1,$ $r_2,$ $0<r_1<r_2< r_0=r(x_0)$ и любых континуумов
$E_1\subset \overline{B(x_0, r_1)}\cap D,$ $E_2\subset
\left(X\setminus B(x_0, r_2)\right)\cap D$ отображение $f$
удовлетворяет соотношению
\begin{equation}\label{eq3*!!A}
 M_p\left(f\left(\Gamma\left(E_1,\,E_2,\,D\right)\right)\right)\ \leqslant
\int\limits_{A} Q(x)\cdot \eta^q(d(x, x_0))\ d\mu(x)
\end{equation}
для каждой измеримой функции $\eta :(r_1,r_2)\rightarrow [0,\infty
],$ такой что имеет место соотношение (\ref{eq*3!!}). Имеет место
следующая теорема.

\medskip
\begin{theorem}\label{th1E} {\sl\, Предположим, $G$ -- область
в метрическом пространстве $(X, d, \mu)$ с локально конечной
борелевской мерой $\mu$ и конечной хаусдорфовой размерностью
$\alpha\geqslant 2,$ а $(X^{\,\prime}, d^{\,\prime},
\mu^{\,\prime})$ -- метрическое пространство, являющееся
$\alpha^{\,\prime}$-регулярным по Альфорсу, в котором выполнено $(1;
p)$-неравенство Пуанкаре, $p\in (\alpha^{\,\prime}-1,
\alpha^{\,\prime}].$ Пусть также область $G$ локально линейно связна
в точках границы, а область $G^{\,\prime}\subset B_R$ является
$QED$-областью относительно $p$-модуля, где $B_R$ -- некоторый шар в
$X^{\,\prime},$ такой что $\overline{B_R}$ -- компакт в
$X^{\,\prime}.$ Тогда произвольный кольцевой $Q$-гомеоморфизм
$f:G\rightarrow G^{\,\prime}$ в точке $b\in
\partial G$ относительно $p$ и $q$-модулей
такой, что $f(G)=G^{\,\prime},$ $q\in (1, \alpha],$ имеет
непрерывное продолжение в точку $b$ при условии, что $Q\in FMO(b).$
}
\end{theorem}

\medskip
Обозначим через $\frak{R}_{p, q, Q, a_0, b_0}(G, G^{\,\prime})$
семейство, состоящее из всех кольцевых $Q$-го\-ме\-о\-мор\-физ\-мов
$f\colon G\rightarrow G^{\,\prime}$ в каждой точке $x_0\in
\overline{G}$ относительно $p$ и $q$-модулей таких, что
$f(a_0)=a_1\ne b_1=f(b_0),$ $f(G)=G^{\,\prime}.$ Основной результат
настоящей работы заключает в себе следующее утверждение.

\medskip
\begin{theorem}\label{th2E}{\sl\,
Пусть $G$ -- область в локально связном и сепарабельном 
метрическом пространстве $(X, d, \mu)$ с конечной хаусдорфовой
размерностью $\alpha\geqslant 2,$ локально линейно связная в каждой
точке своей границы, а $(X^{\,\prime}, d^{\,\prime},
\mu^{\,\prime})$ -- метрическое пространство,
$\alpha^{\,\prime}$-регулярное по Альфорсу, в котором выполнено $(1;
p)$-неравенство Пуанкаре, $p\in (\alpha^{\,\prime}-1,
\alpha^{\,\prime}].$ Предположим, область $G^{\,\prime}\subset B_R$
является $QED$-областью относительно $p$-модуля,
$G^{\,\prime}\subset B_R,$ $B_R$ -- некоторый фиксированный шар в
$X^{\,\prime},$ такой, что $\overline{B_R}$ компакт, причём найдётся
невырожденный континуум $K\subset B_R\setminus G^{\,\prime}.$

Пусть $Q\in FMO(\overline{G})$ и $q\in (1, \alpha].$ Тогда каждое
отображение семейства $\frak{R}_{p, q, Q, a_0, b_0}(G,
G^{\,\prime})$ продолжается по непрерывности на $\partial G,$ при
этом, семейство $\frak{R}_{p, q, Q, a_0, b_0}(\overline{G},
\overline{G^{\,\prime}}),$ состоящее из всех продолженных таким
образом отображений $f:\overline{G}\rightarrow
\overline{G^{\,\prime}},$ является равностепенно непрерывным в
каждой точке $x_0\in \overline{G}.$}
 \end{theorem}

\medskip
{\bf 2. О равностепенной непрерывности го\-ме\-о\-мо\-р\-фи\-з\-мов
внутри области.} Результаты настоящего раздела установлены в работе
\cite{Sev$_4$} в частном случае $p=\alpha,$ $q=\alpha^{\,\prime},$
где $\alpha, \alpha^{\,\prime}$ -- хаусдорфовы размерности
пространств $X$ и $X^{\,\prime},$ соответственно. Как известно,
эффективным методом исследования кольцевых $Q$-отображений в ${\Bbb
R}^n$ является метод <<сингулярных параметров>>, т.е., метод,
позволяющий связать поведение заданной характеристики $Q(x)$ с
поведением модуля соответствующего семейства кривых. Применим этот
метод и к исследованию отображений в метрических пространствах.
Следующая лемма может быть полезной при исследовании свойства
равностепенной непрерывности $Q$-го\-ме\-о\-мор\-физ\-мов,
удовлетворяющих (\ref{eq1}), в наиболее общей ситуации.

\medskip
\begin{lemma}\label{lem4}{\sl\,
Пусть $G$ -- область в метрическом пространстве $(X, d, \mu)$ с
конечной хаусдорфовой размерностью $\alpha\geqslant 2,$ а
$(X^{\,\prime}, d^{\,\prime}, \mu^{\,\prime})$ -- также, некоторое
метрическое пространство с конечной хаусдорфовой размерностью
$\alpha^{\,\prime}\geqslant 2.$ Пусть $p, q \geqslant 1$ и $f\colon
G\rightarrow X^{\,\prime}$ --- кольцевой $Q$-гомеоморфизм в точке
$x_0\in G$ относительно $p$ и $q$-модулей, отображающий область $G$
на некоторую область $f(G)\subset X^{\,\prime},$ кроме того, пусть
$r_0>0$ таково, что шар $B(x_0, r_0)$ лежит со своим замыканием в
$G.$ Предположим, что для некоторого числа $0<\varepsilon_0<r_0,$
некоторого $\varepsilon_0^{\,\prime}\in (0, \varepsilon_0)$ и
семейства неотрицательных измеримых по Лебегу функций
$\{\psi_{\varepsilon}(t)\},$ $\psi_{\varepsilon}\colon (\varepsilon,
\varepsilon_0)\rightarrow [0, \infty],$ $\varepsilon\in\left(0,
\varepsilon_0^{\,\prime}\right),$ выполнено условие
 \begin{equation} \label{eq3.7B}
\int\limits_{\varepsilon<d(x, x_0)<\varepsilon_0}
Q(x)\cdot\psi_{\varepsilon}^q(d(x, x_0))\, d\mu(x)\leqslant
F(\varepsilon, \varepsilon_0)\qquad\forall\,\,\varepsilon\in(0,
\varepsilon_0^{\,\prime}),
 \end{equation}
где $F(\varepsilon, \varepsilon_0)$ --- некоторая заданная функция и
\begin{equation}\label{eq3AB} 0<I(\varepsilon, \varepsilon_0):=
\int\limits_{\varepsilon}^{\varepsilon_0}\psi_{\varepsilon}(t)dt <
\infty\qquad\forall\,\,\varepsilon\in(0,
\varepsilon_0^{\,\prime}).\end{equation}
Тогда для сфер $S_1=S(x_0, \varepsilon)$ и $S_2=S(x_0,
\varepsilon_0),$ $0<\varepsilon<\varepsilon^{\,\prime}_0$ выполнено
неравенство
\begin{equation}\label{eq3B}
M_p(f(\Gamma(S_1, S_2, A)))\leqslant F(\varepsilon,\varepsilon_0)/
I^{\,q}(\varepsilon, \varepsilon_0)\qquad
\forall\,\,\varepsilon\in\left(0,\,\varepsilon_0^{\,\prime}\right)\,,
\end{equation}
где $A=\{x\in G: \varepsilon<d(x, x_0)<\varepsilon_0\}.$}
\end{lemma}

\begin{proof}
Рассмотрим семейство измеримых функций
$$\eta_{\varepsilon}(t)=\psi_{\varepsilon}(t)/I(\varepsilon,
\varepsilon_0 ), \quad t\in(\varepsilon,\, \varepsilon_0)\,.$$
Заметим, что для $\varepsilon\in (0, \varepsilon_0^{\,\prime})$
выполнено равенство
$\int\limits_{\varepsilon}^{\varepsilon_0}\eta_{\varepsilon}(t)\,dt=1.$
Тогда из определения кольцевого $Q$-гомеоморфизма в точке $x_0$
относительно $p$ и $q$-модулей, и соотношений (\ref{eq3.7B})
получаем неравенство (\ref{eq3B}).~$\Box$
 \end{proof}

\medskip
Справедливо следующее утверждение (см.~\cite[предложение~4.7]{AS}).

\medskip
\begin{proposition}\label{pr2}
{\sl Пусть $X$ --- $\alpha$-регулярное по Альфорсу метрическое
пространство с мерой, в котором выполняется $(1; p)$-неравенство
Пуанкаре, $\alpha\geqslant 1,$ $p\in(\alpha-1, \alpha].$ Тогда для
произвольных континуумов $E$ и $F,$ содержащихся в шаре $B(x_0, R),$
и некоторой постоянной $C>0$ выполняется неравенство
$$M_p(\Gamma(E, F, X))\geqslant \frac{1}{C}\cdot\frac{\min\{{\rm diam}\,E, {\rm diam}\,F\}}{R^{1+p-\alpha}}\,.$$ }
\end{proposition}

Теперь сформулируем и докажем утверждение о равностепенной
непрерывности кольцевых $Q$-гомеоморфизмов между метрическими
пространствами в <<максимальной>> степени общности.

\medskip
\begin{lemma}\label{lem1}{\sl\,
Пусть $G$ -- область в локально связном и локально компактном
метрическом пространстве $(X, d, \mu)$ с конечной хаусдорфовой
размерностью $\alpha\geqslant 2,$ а $(X^{\,\prime}, d^{\,\prime},
\mu^{\,\prime})$ -- метрическое пространство,
$\alpha^{\,\prime}$-регулярное по Альфорсу, в котором выполнено $(1;
p)$-неравенство Пуанкаре, $p\in (\alpha^{\,\prime}-1,
\alpha^{\,\prime}].$

Пусть $r_0>0$ таково, что шар $B(x_0, r_0)$ лежит со своим
замыканием в $G$ и $0<\varepsilon_0<r_0.$ Предположим также, что для
некоторого числа $\varepsilon_0^{\,\prime}\in (0, \varepsilon_0)$ и
семейства неотрицательных измеримых по Лебегу функций
$\{\psi_{\varepsilon}(t)\},$ $\psi_{\varepsilon}\colon (\varepsilon,
\varepsilon_0)\rightarrow (0, \infty),$ $\varepsilon\in\left(0,
\varepsilon_0^{\,\prime}\right),$ выполнено условие~\eqref{eq3.7B},
где некоторая заданная функция $F(\varepsilon, \varepsilon_0)$
удовлетворяет условию $F(\varepsilon,
\varepsilon_0)=o(I^{q}(\varepsilon, \varepsilon_0)),$ а
$I(\varepsilon, \varepsilon_0)$ определяется
соотношением~\eqref{eq3AB}.

Пусть $B_R\subset X^{\,\prime}$
--- некоторый фиксированный шар радиуса $R.$ Обозначим
через $\frak{R}_{p, q, x_0, Q, B_R, \delta}(G)$ семейство кольцевых
$Q$-го\-ме\-о\-мор\-физ\-мов $f\colon G\rightarrow B_R\setminus K_f$
в точке $x_0\in G$ относительно $p$ и $q$-модулей таких, что $q\in
(1, \alpha]$ и $\sup\limits_{x, y\in K_f} d^{\,\prime}(x,
y)\geqslant \delta>0,$ где $K_f\subset B_R$ -- некоторый континуум.
Тогда семейство отображений $\frak{R}_{p, q, x_0, Q, B_R,
\delta}(G)$ является равностепенно непрерывным в точке $x_0\in G.$}
 \end{lemma}

\medskip
\begin{proof}
Пусть $x_0\in G,$ $f\in\frak{R}_{p, q, x_0, Q, B_R, \delta}(G).$
Поскольку пространство $X$ локально связно и локально компактно,
можно выбрать последовательность шаров $B(x_0, \varepsilon_k),$
$k=0,1,2,\ldots,$ $\varepsilon_k\rightarrow 0$ при
$k\rightarrow\infty,$ таких что $V_{k+1}\subset \overline{B(x_0,
\varepsilon_k)}\subset V_k,$ где $V_k$ -- континуумы в $G.$ Заметим,
что $f(V_k)$ и $K_f$ -- континуумы в $B_R$ (в частности, $f(V_k)$ --
континуум как непрерывный образ континуума, см. \cite[теорема~1,
III, \S\, 41 и теорема~3, I, \S\, 46]{Ku}). Тогда по предложению
\ref{pr2} при некоторой постоянной $C>0$ получим:
\begin{equation}\label{eq2}
M_p(K_f, f(V_k), X^{\,\prime}))\geqslant
\frac{1}{C}\cdot\frac{\min\{{\rm diam}\,K_f, {\rm diam}\,
f(V_k)\}}{R^{1+p-\alpha^{\,\prime}}}\,.
 \end{equation}
Для кривой $\gamma:[a, b]\rightarrow X$ положим, как обычно,
$$|\gamma|=\{x\in X: \exists\,t\in [a, b]: \gamma(t)=x\}\,.$$
Заметим, что при $k\geqslant 1$ произвольная кривая $\gamma\in
\Gamma(K_f, f(V_k), X^{\,\prime})$ соединяет $f(B(x_0,
\varepsilon_0))$ и $X^{\,\prime}\setminus f(B(x_0, \varepsilon_0)),$
поэтому найдётся точка $y_1\in |\gamma|\cap f(S(x_0,
\varepsilon_0))$ и $t_1\in (0, 1)$ такие, что $\gamma(t_1)=y_1$ и
$|\gamma|_{[0, t_1)}|\subset f(B(x_0, \varepsilon_0))$ (см.
\cite[предложение~13.3]{MRSY} либо \cite[теорема~1.I.46]{Ku}).
Обозначим $\gamma_1:=\gamma|_{[0, t_1)},$ и пусть
$\alpha_1=f^{\,-1}(\gamma_1).$ Заметим, что $|\alpha_1|\subset
B(x_0, \varepsilon_0).$ Заметим далее, что $\alpha_1$ целиком не
лежит ни в $\overline{B(x_0, \varepsilon_{k-1})},$ ни в
$X\setminus\overline{B(x_0, \varepsilon_{k-1})},$ поэтому найдётся
$t_2\in (0, t_1)$ такое, что $\alpha_1(t_2)\in S(x_0,
\varepsilon_{k-1})$ (см. \cite[теорема~1.I.46]{Ku}) и
$|\alpha_{[t_2, t_1]}|\subset X\setminus\overline{B(x_0,
\varepsilon_{k-1})}.$ Положим $\alpha_2=\alpha_1|_{[t_2, t_1]}.$
Заметим, что $\gamma_2:=f(\alpha_2)$ является подкривой $\gamma.$
Исходя из сказанного,
$$\Gamma(K_f, f(V_k),
X^{\,\prime})>\Gamma(f(S(x_0, \varepsilon_{k-1})), f(S(x_0,
\varepsilon_0)), f(A))\,,$$
где $A=\{x\in X: \varepsilon_{k-1}<d(x, x_0)<\varepsilon_0\},$
откуда ввиду соотношения (\ref{eq32*A})
$$M_p(\Gamma(K_f, f(V_k), X^{\,\prime}))\leqslant $$
\begin{equation}\label{eq3}
\leqslant M_p(\Gamma(f(S(x_0, \varepsilon_{k-1})), f(S(x_0,
\varepsilon_0)), f(A)))\,.
\end{equation}
Тогда из соотношений (\ref{eq2}) и (\ref{eq3}) вытекает, что
$$M_p(\Gamma(f(S(x_0, \varepsilon_{k-1})), f(S(x_0, \varepsilon_0)),
f(A)))\geqslant $$
\begin{equation}\label{eq4}
\geqslant\frac{1}{C}\cdot\frac{\min\{{\rm diam}\,K_f, {\rm diam}\,
f(V_k)\}}{R^{1+p-\alpha^{\,\prime}}}\,.
\end{equation}
С другой стороны, из леммы \ref{lem4} и условия $F(\varepsilon,
\varepsilon_0)=o(I^q(\varepsilon, \varepsilon_0))$ вытекает, что
$$M_p(\Gamma(f(S(x_0, \varepsilon_{k-1})), f(S(x_0, \varepsilon_0)),
f(A)))\rightarrow 0$$
при $k\rightarrow \infty,$ поэтому для любого $\sigma>0$ найдётся
$k_0\in {\Bbb N}=k_0(\sigma)$ такой, что при всех $k\geqslant k_0$
$$M_p(\Gamma(f(S(x_0, \varepsilon_{k-1})), f(S(x_0, \varepsilon_0)),
f(A)))<\sigma\,.$$
В таком случае, из (\ref{eq4}) вытекает, что при указанных $k\in
{\Bbb N}$
\begin{equation}\label{eq5}
\min\{{\rm diam}\,K_f, {\rm diam}\, f(V_k)\}< \sigma.
\end{equation}
Поскольку по условию ${\rm diam}\,K_f\geqslant \delta>0$ для всех
$f$ из рассматриваемого семейства отображений, то $\min\{{\rm
diam}\,K_f, {\rm diam}\, f(V_k)\}={\rm diam}\, f(V_k)$ начиная с
некоторого номера $k_1\geqslant k_0,$ $k_1=k_1(\sigma).$ Тогда из
(\ref{eq5}) вытекает, что
\begin{equation}\label{eq1C}
{\rm diam}\, f(V_k)<\sigma
\end{equation}
при всех $k\geqslant k_1.$ Из включений $V_{k+1}\subset
\overline{B(x_0, \varepsilon_k)}\subset V_k$ следует, что
неравенство (\ref{eq1C}) выполнено также в шаре $\overline{B(x_0,
\varepsilon_k)}$ при $k\geqslant k_1(\sigma).$ Положим
$\varepsilon(\sigma):=\varepsilon_{k_1}.$ Окончательно, для числа
$\sigma>0$ найдётся $\varepsilon(\sigma)>0$ такое, что при $d(x,
x_0)< \varepsilon(\sigma)$ выполнено $d^{\,\prime}(f(x),
f(x_0))<\sigma,$ что и означает равностепенную непрерывность
семейства отображений $\frak{R}_{p, q, x_0, Q, B_R, \delta}(G)$ в
точке $x_0.$~$\Box$
\end{proof}

\medskip
{\it Доказательство теоремы~{\em\ref{theor4*!}}} вытекает из
леммы~{\em\ref{lem1}} и \cite[лемма~13.2]{MRSY}.

\medskip
Действительно, согласно \cite[лемма~13.2]{MRSY}, условие $Q\in
FMO(y_0)$ влечёт, что при некотором $\varepsilon_0>0$ и
$\varepsilon\rightarrow 0$
$$\int\limits_{\varepsilon<d(y, y_0)<\varepsilon_0}
Q(x)\cdot\psi^{q}(d(y, y_0))\,
d\mu(y)=O\left(\log\log\frac{1}{\varepsilon}\right)\,,$$
где
$\psi(t)=\left(t\,\log{\frac1t}\right)^{-\alpha/q}>0$ и
$1<q\leqslant\alpha.$
Как и прежде, определим $I(\varepsilon,
\varepsilon_0):=\int\limits_{\varepsilon}^{\varepsilon_0}\psi(t)dt,$
тогда
$$
I(\varepsilon,
\varepsilon_0)=\int\limits_{\varepsilon}^{\varepsilon_0}\psi(t) dt
>\log{\frac{\log{\frac{1}
{\varepsilon}}}{\log{\frac{1}{\varepsilon_0}}}}.
$$
В таком случае, заключаем, что условия
(\ref{eq3.7B})--(\ref{eq3AB}), фигурирующие в лемме \ref{lem1},
выполнены и, значит, из этой леммы вытекает требуемое
утверждение.~$\Box$

\medskip
{\bf 3. О непрерывном продолжении гомеоморфизмов в метрическом
пространстве.} Аналог следующей леммы доказывался В.И. Рязановым и
Р.Р. Салимовым в работе \cite{RSa} для случая, когда границы
отображённых областей сильно достижимы, либо являются слабо плоскими
(см. также статью \cite{Sm}). Указанные условия на границы мы
заменяем ниже требованием вида (\ref{eq4***}), при этом, здесь
присутствуют также некоторые дополнительные ограничения на сами
метрические пространства. Приведённое ниже утверждение установлено в
\cite{Sm} (см. также \cite{RSa}) в частном случае, когда
$p=\alpha^{\,\prime},$ $q=\alpha.$ В случае произвольных $p$ и $q$
наличие упомянутой связи, насколько нам известно, не установлено.

\medskip
\begin{lemma}\label{lem1E} {\sl\, Пусть $G$ -- область 
в метрическом пространстве $(X, d, \mu)$ с конечной хаусдорфовой
размерностью $\alpha\geqslant 2,$ а $(X^{\,\prime}, d^{\,\prime},
\mu^{\,\prime})$ -- метрическое пространство, являющееся
$\alpha^{\,\prime}$-регулярным по Альфорсу, в котором выполнено $(1;
p)$-неравенство Пуанкаре, $p\in (\alpha^{\,\prime}-1,
\alpha^{\,\prime}].$ Пусть также область $G$ локально линейно связна
в точках границы, а область $G^{\,\prime}\subset B_R$ является
$QED$-областью относительно $p$-модуля, где $B_R$ -- некоторый шар в
$X^{\,\prime},$ такой что $\overline{B_R}$ -- компакт в
$X^{\,\prime}.$

Предположим также, что найдётся $\varepsilon_0>0$ и некоторая
положительная измеримая функция $\psi(t),$ $\psi:(0,
\varepsilon_0)\rightarrow (0,\infty),$ такая что для всех
$\varepsilon\in(0, \varepsilon_0)$
\begin{equation}\label{eq7***}
0<I(\varepsilon,
\varepsilon_0)=\int\limits_{\varepsilon}^{\varepsilon_0}\psi(t)dt <
\infty
\end{equation}
и при $\varepsilon\rightarrow 0$ и некотором $q\in (1, \alpha]$
\begin{equation}\label{eq5***}
\int\limits_{A(b, \varepsilon, \varepsilon_0)}
Q(x)\cdot\psi^{\,q}(d(x, b))
 \ d\mu(x) =o(I^{q}(\varepsilon, \varepsilon_0))\,,
\end{equation}
где $A:=A(b, \varepsilon, \varepsilon_0)$ определено в
(\ref{eq49***}). Тогда произвольное кольцевое $Q$-отображение
$f:G\rightarrow G^{\,\prime}$ в точке $b\in
\partial G$ относительно $p$ и $q$-модулей, такое, что $f(G)=G^{\,\prime},$ имеет непрерывное продолжение в точку
$b.$ }
\end{lemma}

\medskip
\begin{proof} Поскольку $G^{\,\prime}\subset B_R$ и $\overline{B_R}$ -- компакт в
$X^{\,\prime},$ предельное множество $C(f, b)$ не пусто.

Предположим противное, а именно, что отображение $f$ не имеет
непрерывного продолжения в точку $b.$ Тогда найдутся, по крайней
мере, две последовательности $x_i,$ $x_i^{\,\prime}\in G,$
$i=1,2,\ldots,$ такие, что $x_i\rightarrow b,$
$x_i^{\,\prime}\rightarrow b$ при $i\rightarrow \infty,$
$f(x_i)\rightarrow y,$ $f(x_i^{\,\prime})\rightarrow y^{\,\prime}$
при $i\rightarrow \infty$ и $y^{\,\prime}\ne y.$ Отметим, что $y$ и
$y^{\,\prime}\in
\partial D^{\,\prime},$ так как $f$ -- гомеоморфизм (см. \cite[предложение~13.5]{MRSY}).
Заметим, что в этом случае найдётся $\delta>0$ такое, что
$d^{\,\prime}(f(x_i), f(x_i^{\,\prime}))\geqslant \delta >0$ при
всех $i\in {\Bbb N}.$ Соединим точки $x_i$ и $x_i^{\,\prime}$ кривой
$C_i,$ целиком лежащей в $B(b, 2^{\,-i}),$ что возможно ввиду
локальной связности области $G$ в точке $b.$ Пусть
$C_i^{\,\prime}=f(C_i),$ тогда ${\rm diam\,}C_i^{\,\prime}\geqslant
\delta>0$ при всех $i\in {\Bbb N}.$ Не ограничивая общности
рассуждений, можно считать, что $G\setminus \overline{B(b,
\varepsilon_0)}\ne \varnothing.$ Выберем произвольную точку $z_0\in
G\setminus \overline{B(b, \varepsilon_0)}$ и соединим её с точкой
$b$ локально спрямляемой кривой, лежащей в $G$ (что возможно ввиду
\cite[предложение~13.2]{MRSY}). Тогда у этой кривой существует
подкривая, лежащая в $G\setminus \overline{B(b, \varepsilon_0)}$
ввиду \cite[предложение~13.3]{MRSY} (см. по этому поводу также
\cite[теорема~1.I.46]{Ku}). Эту подкривую обозначим через $K$ и
заметим, что она представляет собой некоторый фиксированный
континуум в $G\setminus \overline{B(b, \varepsilon_0)};$ тогда при
больших $i\in {\Bbb N}$ имеем: $K\subset G\setminus \overline{B(b,
2^{\,-i})}.$

\medskip
Тогда, с одной стороны, поскольку $G^{\,\prime}$ является
$QED$-областью относительно $p$-модуля, то для некоторой постоянной
$A\geqslant 1$
\begin{equation}\label{eq1E}
M_p(\Gamma(C^{\,\prime}_i, f(K), G^{\,\prime}))\geqslant
\frac{1}{A}\cdot M_p(\Gamma(C^{\,\prime}_i, f(K), X^{\,\prime})\,.
\end{equation}
Поскольку $X^{\,\prime}$ является $\alpha^{\,\prime}$-регулярным по
Альфорсу и, кроме того, в $X^{\,\prime}$ выполнено $(1;
p)$-неравенство Пуанкаре, ввиду предложения \ref{pr2}
\begin{equation}\label{eq2E}
M_p(\Gamma(C^{\,\prime}_i, f(K), X^{\,\prime})\geqslant
\frac{1}{C}\cdot\frac{\min\{{\rm diam}\,f(K), {\rm
diam}\,C^{\,\prime}_i\}}{R^{1+p-\alpha^{\,\prime}}}\geqslant
\delta_1>0\,,
\end{equation}
где $\delta_1$ не зависит от $i.$ Из (\ref{eq1E}) и (\ref{eq2E})
вытекает, что
\begin{equation}\label{eq3E}
M_p(\Gamma(C^{\,\prime}_i, f(K),
G^{\,\prime}))\geqslant\delta_2>0\,,
\end{equation}
где $\delta_2$ не зависит от $i.$

\medskip
С другой стороны, рассмотрим семейство кривых $\Gamma_i,$
соединяющих $K$ и $C_i.$ Рассмотрим функцию
$$\eta(t)=\left\{
\begin{array}{rr}
\psi(t)/I(2^{-i}, \varepsilon_0), &   t\in (2^{-i},
\varepsilon_0),\\
0,  &  t\in {\Bbb R}\setminus (2^{-i}, \varepsilon_0)\,,
\end{array}
\right. $$ где $I(\varepsilon,
\varepsilon_0):=\int\limits_{\varepsilon}^{\varepsilon_0}\psi(t)dt,$
удовлетворяет условию нормировки вида (\ref{eq*3!!}) при
$r_1:=2^{-i},$ $r_2:=\varepsilon_0.$ В силу определения кольцевого
$Q$-отображения $f:G\rightarrow G^{\,\prime}$ в точке $b\in
\partial G$ относительно $p$ и $q$-модулей, а также соотношений (\ref{eq7***})--(\ref{eq5***}), будем
иметь:
\begin{equation}\label{eq11*}
M_p\left(f\left(\Gamma_i\right)\right)=M_p(\Gamma(C^{\,\prime}_i,
f(K), G^{\,\prime}))\leqslant \Delta(i)\,,
\end{equation}
где $\Delta(i)\rightarrow 0$ при $i\rightarrow \infty.$ Однако,
(\ref{eq11*}) противоречит (\ref{eq3E}), что и доказывает
утверждение леммы.~$\Box$\end{proof}

\medskip
{\it Доказательство теоремы \ref{th1E}} вытекает из леммы
\ref{lem1E} на основании рассуждений, аналогичных рассуждениям,
сделанных при доказательстве теоремы \ref{theor4*!}.~$\Box$

\medskip
{\bf 4. О равностепенной непрерывности го\-ме\-о\-мо\-р\-физ\-мов в
замыкании области.} Обозначим через $\frak{R}_{p, q, Q, a_0, b_0}(G,
G^{\,\prime})$ семейство, состоящее из всех ко\-ль\-цевых $
Q$-гомеоморфизмов $f\colon G\rightarrow G^{\,\prime}$ в каждой точке
$x_0\in \overline{G}$ относительно $p$ и $q$-модулей, таких, что
$f(a_0)=a_1\ne b_1=f(b_0),$ $f(G)=G^{\,\prime}.$  Имеет место
следующее утверждение.

\medskip
\begin{lemma}\label{lem2E}{\sl\, Пусть $G$ -- область в локально связном и сепарабельном 
метрическом пространстве $(X, d, \mu)$ с конечной хаусдорфовой
размерностью $\alpha\geqslant 2,$ локально линейно связная в каждой
точке своей границы, а $(X^{\,\prime}, d^{\,\prime},
\mu^{\,\prime})$ -- метрическое пространство,
$\alpha^{\,\prime}$-регулярное по Альфорсу, в котором выполнено $(1;
p)$-неравенство Пуанкаре, $p\in (\alpha^{\,\prime}-1,
\alpha^{\,\prime}].$ Предположим, область $G^{\,\prime}\subset B_R$
является $QED$-областью относительно $p$-модуля, $B_R$ -- некоторый
фиксированный шар в $X^{\,\prime},$ $\overline{B_R}$ -- компакт в
$X^{\,\prime},$ причём найдётся невырожденный континуум $K\subset
B_R\setminus G^{\,\prime}.$

Предположим, для каждой точки $x_0\in \overline{G},$ некоторого
числа $\varepsilon_0^{\,\prime}\in (0, \varepsilon_0),$
$\varepsilon_0=\varepsilon_0(x_0),$ и семейства
$\{\psi_{\varepsilon}(t)\}$ измеримых по Лебегу функций
$\psi_{\varepsilon}\colon (\varepsilon, \varepsilon_0)\rightarrow
(0, \infty),$ $\varepsilon\in\left(0,
\varepsilon_0^{\,\prime}\right),$ выполнено условие~\eqref{eq3.7B},
где некоторая заданная функция $F(\varepsilon, \varepsilon_0)$
удовлетворяет условию $$F(\varepsilon,
\varepsilon_0)=o(I^{q}(\varepsilon, \varepsilon_0))\,, q\in (1,
\alpha]\,,$$ а $I(\varepsilon, \varepsilon_0)$ определяется
соотношением~\eqref{eq3AB}.

Тогда каждое отображение семейства $\frak{R}_{p, q, Q, a_0, b_0}(G,
G^{\,\prime})$ продолжается по непрерывности на $\partial G,$ при
этом, семейство отображений $\frak{R}_{p, q, Q, a_0,
b_0}(\overline{G}, \overline{G^{\,\prime}}),$ состоящее из всех
продолженных таким образом отображений $f:\overline{G}\rightarrow
\overline{G^{\,\prime}}$ является равностепенно непрерывным в каждой
точке $x_0\in \overline{G}.$
}
 \end{lemma}

\medskip
\begin{proof} Заметим, что $G$ -- локально компактно. Действительно,
так как $\overline{G^{\,\prime}}\subset B_R$ и $\overline{B_R}$ --
компакт в $X^{\,\prime},$ то $\overline{G^{\,\prime}}$ также компакт
как замкнутое подмножество компактного пространства
$\overline{B_R}.$ Тогда для достаточно малого $\varepsilon>0$ шары
$\overline{B(y_0, \varepsilon)}$ компактны при $y_0\in
\overline{G^{\,\prime}}.$ С другой стороны, поскольку
$f(G)=G^{\,\prime}$ при произвольном $f\in \frak{R}_{p, q, Q, a_0,
b_0}(G, G^{\,\prime}),$ то какова бы ни была точка $x_0\in G,$ для
точки $y_0=f(x_0)$ множество $f^{\,-1}(\overline{B(y_0,
\varepsilon)})$ компактно как непрерывный образ компакта и,
одновременно, является окрестностью точки $x_0.$

Кроме того, заметим, что $G^{\,\prime}$ имеет хаусдорфову
размерность $\alpha^{\,\prime},$ что вытекает из
$\alpha^{\,\prime}$-регулярности по Альфорсу пространства
$X^{\,\prime}$ (см.\ рассуждения на с.~61 в~\cite{He}).

В таком случае, равностепенная непрерывность семейства продолженных
отображений $\frak{R}_{p, q, Q, a_0, b_0}(\overline{G},
\overline{G^{\,\prime}})$ во внутренних точках области $G$ является
утверждением леммы \ref{lem1}, а возможность непрерывного на границу
по непрерывности -- леммы \ref{lem1E}. Осталось доказать
равностепенную непрерывность семейства $\frak{R}_{p, q, Q, a_0,
b_0}(\overline{G}, \overline{G^{\,\prime}})$ в точках $\partial G.$

Предположим противное, тогда найдётся $x_0\in
\partial G$ и число $a>0$ такое, что для каждого
$m=1,2,\ldots$ существуют точка $x_m\in \overline{G}$ и элемент
$f_m$ семейства $\frak{R}_{p, q, Q, a_0, b_0}(\overline{G},
\overline{G^{\,\prime}})$ такие, что $d(x_0, x_m)< 1/m$ и
\begin{equation}\label{eq6***}
d^{\,\prime}\left(f_m(x_m), f_m(x_0)\right)\geqslant a\,.
\end{equation}
%
Ввиду возможности непрерывного продолжения каждого $f_m$ на границу
$G,$ мы можем считать, что $x_m\in G.$

В силу локальной линейной связности области $G$ в точке $x_0$
найдётся последовательность окрестностей $V_m$ точки $x_0$ с ${\rm
diam}\,V_m\rightarrow 0$ при $m\rightarrow\infty,$ такие что
множества $G\cap V_m $ являются областями и $x_m\in G\cap V_m.$ Т.к.
граничные точки области, локально связной на границе являются
достижимыми из $G$ некоторым локально спрямляемым путём, см.
\cite[предложение~13.2]{MRSY}, мы можем соединить точки $x_m$ и
$x_0$ кривой $\gamma_m(t):[0,1]\rightarrow G$ такой, что
$\gamma_m(0)=x_0,$ $\gamma_m(1)=x_m$ и $\gamma_m(t)\in V_m$ при
$t\in (0,1).$ Можно считать, что $|\gamma_m|\subset B(x_0, 1/m).$
Обозначим через $C_m$ образ кривой $\gamma_m(t)$ при отображении
$f_m.$ Из соотношения (\ref{eq6***}) вытекает, что
\begin{equation}\label{eq5.1}
{\rm diam}\,C_m\geqslant a\qquad\forall\, m\in {\Bbb N}\,.
\end{equation}

\medskip
Соединим точки $a_0$ и $b_0\in G$ из условия леммы кривой $\beta:[0,
1]\rightarrow G,$ лежащей в $G,$ такой что $\beta(0)=a_0$ и
$\beta(1)=b_0.$ Можно считать, что ${\rm dist}\, (|\beta|,
\partial G)>\varepsilon_0(x_0),$ где $\varepsilon_0(x_0)$ -- из условия леммы,
а, как обычно, $|\beta|=\{x\in G: \exists\,t\in [0, 1]:
\beta(t)=x\}$ -- носитель кривой $\beta.$ Поскольку семейство
$\frak{R}_{p, q, Q, a_0, b_0}(\overline{G},
\overline{G^{\,\prime}})$ является равностепенно непрерывным в $G,$
$\overline{B_R}$ -- компакт, а $G$ -- сепарабельно, $\frak{R}_{p, q,
Q, a_0, b_0}(\overline{G}, \overline{G^{\,\prime}})$ является
нормальным ввиду критерия Арцела--Асколи (см.
\cite[пункт~20.4]{Va}). В таком случае, не ограничивая общности
рассуждений, можно считать, что найдётся непрерывное отображение
$f:G\rightarrow B_R,$ такое что $\sup\limits_{x\in
|\beta|}d^{\,\prime}(f_m(x), f(x))\rightarrow 0$ при $m\rightarrow
\infty.$ Тогда $f(|\beta|)$ -- компакт в $X^{\,\prime}$ как образ
компакта $|\beta|\subset G$ при непрерывном отображении $f.$

\medskip
Возможны две ситуации: 1) $f(|\beta|)\subset G^{\,\prime},$ тогда
полагаем $B:=f(|\beta|);$ 2) $f(|\beta|)\cap \partial
G^{\,\prime}\ne \varnothing.$ В этом случае полагаем
$t_0:=\{\sup\limits_{t\in [0, 1]}t:$ $f(\beta(r))\in G^{\,\prime}$
при всех $r\in[0, t]\}.$ Возьмём теперь произвольное $s_0<t_0$ и
положим $B:=f(|\beta_{[0, s_0]}|).$ Очевидно, в обеих из двух
ситуаций $B$ -- невырожденный континуум в $G^{\,\prime},$ при этом,
существует компакт $C=|\beta_{[0, s_0]}|,$ $0<s_0\leqslant 1,$ такой
что $f(C)=B.$ Тогда ввиду локально равномерной сходимости $f_m$ к
$f$ найдётся $k_0\in {\Bbb N}$ такой, что
\begin{equation}\label{eq5E}
{\rm diam}\,(f_m(C))>\frac{{\rm
diam}\,(f(C))}{2}:=\delta_0>0\quad\forall\,m\geqslant k_0\,.
\end{equation}
Рассмотрим теперь семейство кривых $\Gamma_m=\Gamma(f_m(C), C_m,
G^{\,\prime}).$ Поскольку по условию $G^{\,\prime}$ является
$QED$-областью, а пространство $X^{\,\prime}$ является
$\alpha^{\,\prime}$-регулярным по Альфорсу, в котором выполнено $(1;
p)$-неравенство Пуанкаре, то ввиду предложения \ref{pr2},
(\ref{eq5.1}) и (\ref{eq5E}), будем иметь:
\begin{equation}\label{eq6A}
M_p(\Gamma_m)\geqslant \frac{1}{A}\cdot M_p(\Gamma(f_m(C), C_m,
X^{\,\prime}))\geqslant
\end{equation}
$$\geqslant\frac{1}{AC}\cdot\frac{\min\{{\rm
diam}\,f_m(C), {\rm
diam}\,C_m\}}{R^{1+p-\alpha^{\,\prime}}}\geqslant M_1>0\,,$$
где $M_1$ -- некоторая постоянная, которая не зависит от $m\in {\Bbb
N}.$
С другой стороны, так как $f_m$ -- гомеоморфизм, то
$\Gamma_m=\Gamma(f_m(C), C_m, G^{\,\prime})$ $=f_m(\Gamma(C,
|\gamma_m|, G),$ причём $C\subset G\setminus \overline{B(x_0,
\varepsilon_0)},$ а $|\gamma_m|\subset B(x_0, 1/m).$ Тогда по
определению кольцевого $Q$-гомеоморфизма в точке $x_0$ относительно
$p$ и $q$-модулей
$$M_p(\Gamma_m)=M_p(\Gamma(f_m(C), C_m,
G^{\,\prime}))=M_p(f_m(\Gamma(C), |\gamma_m|, G))\leqslant$$
\begin{equation}\label{eq10}
\leqslant \int\limits_{A(x_0, \frac{1}{m}, \varepsilon_0)} Q(x)\cdot
\eta^q(d(x, x_0)) d\mu(x)
\end{equation}
для каждой измеримой функции $\eta: (\frac{1}{m},
\varepsilon_0)\rightarrow [0,\infty ],$ такой что
$\int\limits_{\frac{1}{m}}^{\varepsilon_0}\eta(r)dr \geqslant 1.$
Заметим, что функция
$$\eta(t)=\left\{
\begin{array}{rr}
\psi(t)/I(1/m, \varepsilon_0), &   t\in (1/m,
\varepsilon_0),\\
0,  &  t\in {\Bbb R}\setminus (1/m, \varepsilon_0)\,,
\end{array}
\right. $$ где $I(\varepsilon,
\varepsilon_0):=\int\limits_{\varepsilon}^{\varepsilon_0}\psi(t)dt,$
удовлетворяет условию нормировки вида (\ref{eq*3!!}) при $r_1:=1/m,$
$r_2:=\varepsilon_0,$ поэтому из условий (\ref{eq10}) и
(\ref{eq3.7B}) вытекает, что
\begin{equation}\label{eq11}
M_p(\Gamma_m)\leqslant \alpha(1/m)\rightarrow 0
\end{equation}
при $m\rightarrow \infty,$ где $\alpha(\varepsilon)$ -- некоторая
неотрицательная функция, стремящаяся к нулю при
$\varepsilon\rightarrow 0,$ которая существует ввиду условия
(\ref{eq3.7B}).
Однако, соотношение (\ref{eq11}) противоречит (\ref{eq6A}).
Полученное противоречие указывает на то, что исходное предположение
(\ref{eq6***}) было неверным, и, значит, семейство отображений
$\frak{R}_{p, q, Q, a_0, b_0}(\overline{G},
\overline{G^{\,\prime}})$ равностепенно непрерывно в точке $x_0\in
\partial G.$~$\Box$
\end{proof}

\medskip
{\it Доказательство теоремы \ref{th2E}} вытекает из леммы
\ref{lem2E} на основании рассуждений, аналогичных рассуждениям,
сделанных при доказательстве теоремы \ref{theor4*!}.~$\Box$

\medskip
{\bf 5. Заключительные замечания.} Большая часть условий,
присутствующих в основных результатах статьи, по-видимому, являются
только достаточными. Тем не менее, опишем некоторые требования,
которые можно обозначить, как <<близкие к необходимым>>.

Прежде всего, в теоремах \ref{theor4*!}, \ref{th1E} и \ref{th2E}
нельзя, вообще говоря, отказаться от условия $Q\in FMO,$ заменив его
более слабым требованием $Q\in L^p,$ каким бы большим ни было число
$p\geqslant 1.$ Для простоты рассмотрим $X=X^{\,\prime}={\Bbb R}^n$
со стандартной евклидовой метрикой и Лебеговой мерой. Пусть, кроме
того, $q=p=n.$ В этом случае отображение, удовлетворяющее оценке
(\ref{eq3*!!A}) будем называть просто <<кольцевым
$Q$-гомеоморфизмом>>. Ниже приведён результат, относящийся к
равностепенной непрерывности (по поводу устранения особенности см.,
напр., \cite[предложение~6.3]{MRSY}).

\medskip
Положим $D:={\Bbb B}^n\setminus \{0\}\subset {\Bbb R}^n,$
$D^{\,\prime}:=B(0, 2)\setminus\{0\} \subset {\Bbb R}^n.$ Обозначим
через $\frak{A}_{Q}$ семейство всех кольцевых $Q$-гомеоморфизмов
$g:{\Bbb B}^n\setminus \{0\}\rightarrow {\Bbb R}^n$ в точке $0.$
Справедливо следующее утверждение.

\medskip
\begin{theorem}\label{th3.10.1}{\sl\, Для каждого $p\geqslant 1$ найдётся
функция $Q:{\Bbb B}^n\rightarrow [1, \infty],$ $Q(x)\in L^p({\Bbb
B}^n)$ и последовательность $g_m\in \frak{A}_{Q}$ такая, что каждый
элемент $g_m$ имеет непрерывное продолжение в точку $x_0=0,$ при
этом, семейство $\left\{g_m(x)\right\}_{m=1}^{\infty}$ не является
равностепенно непрерывным в точке $x_0=0.$}
\end{theorem}

\medskip
{\it Доказательство теоремы \ref{th3.10.1}} может быть найдено в
\cite[теорема~8]{Sev$_4$}.~$\Box$

\medskip
В условиях теоремы \ref{th2E}, даже в случае $Q(x)\equiv 1,$ от
условия фиксации, по крайней мере, одной внутренней точки области
$D$ каждым гомеоморфизмом $f$ соответствующего семейства
отображений, отказаться нельзя. Сказанное выше показывает следующий
пример семейства конформных ($Q(x)\equiv 1$) отображений на
плоскости: $f_t(z)=\frac{z-t}{1-tz},$ которое при каждом
фиксированном $t\in (-1,1)$ переводит область $D={\Bbb B}^2\subset
{\Bbb C}$ на $D^{\,\prime}={\Bbb B}^2\subset {\Bbb C},$ см., напр.,
\cite[соотношение~(12), гл.~V, $\S\,1$]{GK}. При этом, при каждом
фиксированном $z\in {\Bbb B}^2,$ $f_t(z)\rightarrow -1$ при
$t\rightarrow 1,$ в то же время, $f_t(1)=1$ при всех $t\in (-1, 1),$
откуда следует, что семейство $f_t(z)$ не является равностепенно
непрерывным в точке $z_0=1.$ Относительно необходимости фиксации
двух и более точек области мы ничего не можем сказать -- это условие
может относиться к методу доказательства и, таким образом, не быть
необходимым.

КОНТАКТНАЯ ИНФОРМАЦИЯ

\medskip
\noindent{{\bf Севостьянов Евгений Александрович}
\\Житомирский государственный университет имени Ивана Франко
\\
10 008 Украина, г.~Житомир, ул.~Большая Бердичевская, д.~40,
\\кафедра математического анализа, моб. тел.~066 959 50 34 \\
e-mail: esevostyanov2009@gmail.com}

\end{document}